\documentclass[11pt]{amsart}

\usepackage[numbers]{natbib}
\usepackage{amssymb}
\usepackage{amsmath}
\usepackage{amsthm}

\numberwithin{equation}{section}

\newtheorem{theorem}{Theorem}
\newtheorem{proposition}{Proposition}
\newtheorem{corollary}{Corollary}
\newtheorem{example}{Example}
\newtheorem{lemma}{Lemma}
\newtheorem{remark}{Remark}

\newcommand{\Z}{\mathcal{Z}}

\newcommand{\supp}{\mathrm{supp}}
\newcommand{\B}{\mathcal{B}}

\newcommand{\A}{\mathcal{A}}
\newcommand{\realnos}{\mathbb{R}}

\begin{document}

\title[Central Limit Theorems for Maps]{Central Limit Theorems for Non-Invertible Measure Preserving
Maps}
\author[M. C. Mackey]{Michael C. Mackey}
\address{Departments of Physiology, Physics \& Mathematics and Centre for
Nonlinear Dynamics, McGill University, 3655 Promenade Sir William
Osler, Montreal, QC, CANADA, H3G 1Y6}%
\email{michael.mackey@mcgill.ca}%
\author[M. Tyran-Kami\'nska]{Marta Tyran-Kami\'nska${}^\dagger$}%
\address{Institute of Mathematics,
University of Silesia, ul. Bankowa 14, 40-007 Katowice, POLAND}%
\email{mtyran@us.edu.pl}%
\date{\today}
\thanks{${}^\dagger$Corresponding author}%

\maketitle

\begin{abstract}
Using the Perron-Frobenius operator we establish a new functional
central limit theorem result for non-invertible measure preserving
maps that are not necessarily ergodic. We apply the result to
asymptotically periodic transformations and give an extensive
specific example using the tent map.

\end{abstract}

\section{Introduction}

This paper is motivated by the question ``How can we produce the
characteristics of a Wiener process (Brownian motion) from a
semi-dynamical system?".   This question is intimately connected
with central limit theorems for non-invertible maps and various
invariance principles. Many results on  central limit theorems and
invariance principles for maps have been proved, see e.g. the
surveys \citet{denker89} and \citet{mackeytyran}. These results
extend back over some decades,
and include the work of  \citet{boyarsky79}, 
\citet{gouezel04b}, \citet{jab83}, \citet{rousseau-egele}, and
\citet{wong79} for the special case of maps of the unit interval.
Martingale approximations, developed by \citet{gordin69}, were used
by \citet{keller80}, \citet{liverani},
\citet{melbourne04}, \citet{melbourne},
and \citet{tyran}, 
 to give more general results.
%
%

Throughout this paper, $(Y,\B,\nu)$ denotes a probability measure
space and $T:Y\to Y$ a non-invertible measure preserving
transformation. Thus $\nu$ is invariant under $T$ {\it i.e.}
$\nu(T^{-1}(A))=\nu(A)$ for all $A\in \B$. The transfer operator
$\mathcal{P}_{T}:L^1(Y,\B,\nu)\to L^1(Y,\B,\nu)$, by definition,
satisfies
\[
\int \mathcal{P}_{T}f(y)g(y)\nu(dy)=\int f(y)g(T(y))\nu(dy)
\]
for all $f\in L^1(Y,\B,\nu)$ and $g\in L^{\infty}(Y,\B,\nu)$.

Let $h\in L^2(Y,\B,\nu)$ with $\int h(y)\nu(dy)=0$. Define the
process $\{w_n(t):t\in[0,1]\}$ by
\begin{equation}\label{e:wnh}
w_n(t)=\frac{1}{\sqrt{n}}\sum_{j=0}^{[nt]-1}h\circ T^j
\;\;\mbox{for}\;\;t\in [0,1],\;n\ge 1
\end{equation}
(the sum from $0$ to $-1$ is set equal to $0$), where $[x]$ denotes
the integer part of $x$.  
For  each $y$, $w_n(\cdot)(y)$ is
an element of the Skorohod space $D[0,1]$ of all functions which are
right continuous and have left-hand limits
equipped with the Skorohod topology. 
\[
\rho_S(\psi,\widetilde{\psi})=\inf_{s\in
\mathcal{S}}\left(\sup_{t\in[0,1]}|\psi(t)-\widetilde{\psi}(s(t))|+\sup_{t\in[0,1]}|t-s(t)|\right),\;\;\psi,\widetilde{\psi}\in
D[0,1],
\]
where $\mathcal{S}$ is the family of strictly increasing, continuous mappings $s$
of $[0,1]$ onto itself such that $s(0)=0$ and $s(1)=1$ \citep[Section
14]{billingsley68}.

Let $\{w(t):t\in[0,1]\}$ be a standard Brownian motion. Throughout the paper the
notation
\begin{equation*}
w_n\to^d \sqrt{\eta}w,
\end{equation*}
where $\eta$ is a random variable independent of the Brownian
process $w$, denotes the weak convergence of the sequence $w_n$
in the Skorohod space $D[0,1]$.



Our main result, which is proved using techniques similar to those in
\citet{peligradutev} and  \citet{peligradutevwu}, is the following:

\begin{theorem} \label{t:CLT2} Let $T$ be a non-invertible measure-preserving transformation
on the probability space $(Y,\B,\nu)$ and let $\mathcal{I}$ be the
$\sigma$-algebra of all $T$-invariant sets. Suppose  $h\in
L^2(Y,\B,\nu)$ with $\int h(y)\nu(dy)=0$ is such that
\begin{equation} \label{concltpo} \sum_{n=1}^\infty
{n^{-\frac{3}{2}}}\biggl\lVert\sum_{k=0}^{n-1}{\mathcal P}_{T}^k
h\biggr\rVert_2<\infty.
\end{equation}
Then
\begin{equation}\label{e:clt2}
w_n\to^d \sqrt{\eta} w,
\end{equation}
where $\eta=E_\nu(\tilde{h}^2|\mathcal{I})$ 
and $\tilde{h}\in L^2(Y,\B,\nu)$ is such that
$\mathcal{P}_{T}\tilde{h}=0$ and
    $$
    \lim_{n\to\infty}\biggl\lVert\frac{1}{\sqrt{n}}\sum_{j=0}^{n-1}(h-\tilde{h})\circ
T^j\biggr\rVert_2 = 0.
    $$
\end{theorem}

Recall that $T$ is {\it ergodic} (with respect to $\nu$) if, for
each $A\in\B$ with $T^{-1}(A)=A$, we have $\nu(A)\in \{0,1\}$. Thus
if $T$ is ergodic then $\mathcal{I}$  is a trivial $\sigma$-algebra,
so $\eta$ in~\eqref{e:clt2} is a constant random variable.
Consequently, Theorem~\ref{t:CLT2}  significantly generalizes
\citet[Theorem 4]{tyran}, where it was assumed that $T$ is ergodic
and  there is $\alpha<1/2$ such that
\[
\biggl\lVert\sum_{k=0}^{n-1}{\mathcal P}_{T}^k
h\biggr\rVert_2=O(n^\alpha)
\]
(We use the notation $b(n)=O(a(n))$ if $\limsup_{n\to\infty}b(n)/a(n)<\infty$).

Usually, in proving central limit theorems for specific examples of
transformations one assumes that the transformation is mixing. For
non-invertible ergodic transformations for which the transfer
operator is quasi-compact on some subspace $F\subset L^2(\nu)$ with
norm $|\cdot|\ge \lVert \cdot\rVert_2$, the central limit theorem
and its functional version was given in~\citet{melbourne04}. Since
quasicompactness implies exponential decay of the $L^2$ norm, our
result applies, thus extending the results of \citet{melbourne04} to
the non-ergodic case. For examples of transformations in which the
decay of the $L^2$ norm is slower than exponential and our results
apply, see \citet{tyran}.

In the case of invertible transformations, non-ergodic versions of
the central limit theorem and its functional generalizations were
studied in \citet{volny87b, volny87, volny89, volny93} using
martingale approximations. In a recent review by
\citet{merlevedeetal06}, the weak invariance principle was studied
for stationary sequences $(X_k)_{k\in\mathbb{Z}}$ which, in
particular, can be described as $X_k=X_0\circ {T}^k$, where $T$ is a
measure preserving invertible transformation on a probability space
and $X_0$ is measurable with respect to a $\sigma$-algebra
$\mathcal{F}_0$ such that $\mathcal{F}_0\subset
T^{-1}(\mathcal{F}_0)$. Choosing a $\sigma$-algebra $\mathcal{F}_0$
for a specific example of invertible transformation is not an easy
task and the requirement that $X_0$ is $\mathcal{F}_0$-measurable
may sometimes be too restrictive (see
\cite{conze01,liverani}). Sometimes, it is possible to reduce
an invertible transformation to a non-invertible one (see
\cite{melbourne,tyran}). Our result in the non-invertible case
extends \citet[Theorem 1.1]{peligradutev}, which is also to be found
in \citet[Theorem 11]{merlevedeetal06}, where a condition
introduced by \citet{maxwell00} is assumed. 
In \citet{tyran} the condition was transformed to Equation \eqref{concltpo}. In
the proof of our result we use Theorem~4.2 in \citet{billingsley68} and
approximation techniques which was motivated by \citet{peligradutev}. The
corresponding maximal inequality in our non-invertible setting is stated in
Proposition~1 and its proof, based on ideas of \citet{peligradutevwu}, is provided
in Appendix~\ref{a:app} for completeness. As in \citet{peligradutev}, the random
variable $\eta$ in Theorem~\ref{t:CLT2} can also be obtained as a limit in $L^1$,
which we state in Appendix~\ref{b:app}.

The outline of the paper is as follows.  Following the presentation of some
background material in Section \ref{s:pre}, we turn to a proof of our main result
Theorem~\ref{t:CLT2} in Section \ref{s:CLT2}. Section \ref{s:aptrans} introduces
asymptotically periodic transformations as a specific example of a system to which
Theorem \ref{t:CLT2} applies.  We analyze the specific example of an
asymptotically periodic family of tent maps  in Section~\ref{s:tent}.

\section{Preliminaries}\label{s:pre}
The definition of the Perron-Frobenius (transfer) operator for $T$
depends on a given $\sigma$-finite measure $\mu$ on the measure
space $(Y,\B)$ with respect to which $T$ is {nonsingular}, {\it
i.e.} $\mu(T^{-1}(A))=0$ for all $A\in\B$ with $\mu(A)=0$. Given
such a measure the {\it transfer operator} $P:L^1(Y,\B,\mu)\to
L^1(Y,\B,\mu)$ is defined as follows. For any $f\in L^1(Y,\B,\mu)$,
there is a unique element $P f$ in $L^1(Y,\B,\mu)$ such that
    \begin{equation}
    \int_{A}Pf(y)\mu(dy)=\int_{T^{-1}(A)}f(y)\mu(dy)\qquad\mbox{for all }A\in\B.
    \label{fpoper}
    \end{equation}
This in turn gives rise to different operators for different
underlying measures on $\B$. Thus if $\nu$ is invariant for $T$,
then $T$ is nonsingular and the transfer operator ${\mathcal
P}_{T}:L^1(Y,\B,\nu)\to L^1(Y,\B,\nu)$ is well defined. Here we
write ${\mathcal P}_{T}$ to emphasize that the underlying measure
$\nu$ is invariant under $T$.

The Koopman operator is defined by
\[
U_Tf=f\circ T
\]
for every measurable $f:Y\to\realnos$. In particular, $U_T$ is also
well defined for $f\in L^1(Y,\B,\nu)$ and is an isometry of
$L^1(Y,\B,\nu)$ into $L^1(Y,\B,\nu)$, {\it i.e.} $||U_T
f||_1=||f||_1$ for all $f\in L^1(Y,\B,\nu)$. Since the measure $\nu$
is finite, we have $L^p(Y,\B,\nu)\subset L^1(Y,\B,\nu)$ for $p\ge
1$. The operator $U_T:L^p(Y,\B,\nu)\to L^p(Y,\B,\nu)$ is also an
isometry on $L^p(Y,\B,\nu)$.



The following relations hold between the operators $U_T,{\mathcal
P}_{T}\colon L^1(Y,\B,\nu)\to L^1(Y,\B,\nu)$
 \begin{equation}
{\mathcal P}_{T} U_T f=f\;\;\mbox{and}\;\;U_T {\mathcal P}_{T}
f=E_{\nu}(f|T^{-1}(\B))\label{fpcond}
 \end{equation} for $f\in L^1(Y,\B,\nu),$ where
$E_{\nu}(\cdot|T^{-1}(\B))\colon L^1(Y,\B,\nu)\to L^1(Y,T^{-1}(\B),\nu)$ is the
operator of conditional expectation. Note that if the transformation $T$ is
invertible then $U_T {\mathcal P}_{T}f=f$ for $f\in L^1(Y,\B,\nu)$.  


\begin{theorem} \label{t:CLT1} Let $T$ be a non-invertible measure-preserving transformation
on the probability space $(Y,\B,\nu)$ and let $\mathcal{I}$ be the
$\sigma$-algebra of all $T$-invariant sets. Suppose that $h\in
L^2(Y,\B,\nu)$ is such that ${\mathcal P}_{T}h=0$. Then
\begin{equation*}
    w_n\to^d
    \sqrt{\eta}w,
\end{equation*}
where $\eta=E_\nu(h^2|\mathcal{I})$ is a random variable independent
of the Brownian motion $\{w(t):t\in[0,1]\}$.
%
%

\end{theorem}
\begin{proof}
When $T$ is ergodic,  a direct proof based on the fact that the
family
\[
\{T^{-n+j}(\B),\frac{1}{\sqrt{n}}h\circ T^{n-j}: 1\le j\le n, n\ge
1\}
\]
is a martingale difference array is given in \citet[Appendix
A]{mackeytyran} and  uses the Martingale Central Limit Theorem (cf.
\citet[Theorem 35.12]{billingsley95}) together with the Birkhoff
Ergodic Theorem.  This can be extended to the case of non-ergodic
$T$  by using a version of the Martingale Central Limit Theorem due
to \citet[Corollary p. 561]{eagleson75}.
\end{proof}

\begin{example}
We illustrate Theorem~\ref{t:CLT1} with an example. Let
$T:[0,1]\to[0,1]$ be defined by
\[
T(y)=\left\{
\begin{array}{ll}
2y , & y \in  [0,\frac{1}{4} )\\ 2y -\frac{1}{2}, & y
\in  [\frac{1}{4}, \frac{3}{4} ),\\
2y-1, & y\in[\frac{3}{4},1].
\end{array}
\right.
\]
Observe that the Lebesgue measure on $([0,1],\B([0,1]))$ is
invariant for $T$ and that $T$ is not ergodic since
$T^{-1}([0,\frac{1}{2}])=[0,\frac{1}{2}]$ and
$T^{-1}([\frac{1}{2},1])=[\frac{1}{2},1]$. The transfer operator is
given by
\[
\mathcal{P}_Tf(y)=\frac{1}{2}f\left(\frac{1}{2}y\right)1_{[0,\frac{1}{2})}(y)
+ \frac{1}{2}f\left(\frac{1}{2}y+\frac{1}{4}\right)+
\frac{1}{2}f\left(\frac{1}{2}y +
\frac{1}{2}\right)1_{[\frac{1}{2},1]}(y).
\]
Consider the function
\[
h(y)=\left\{
\begin{array}{ll}
 1 , & y \in  [0,\frac{1}{4} )\\ - 1, & y
\in  [\frac{1}{4}, \frac{1}{2} ),\\
- 2, & y
\in  [\frac{1}{2}, \frac{3}{4} ),\\
 2, & y\in[\frac{3}{4},1].
\end{array}
\right.
\]
A straightforward calculation shows that $\mathcal{P}_Th=0$ and
 $E_\nu(h^2|\mathcal{I})=1_{[0,\frac{1}{2}]}+4 \,1_{[\frac{1}{2},1]}$. Thus
Theorem~\ref{t:CLT1} shows that
\[
w_n\to^d \sqrt{E_\nu(h^2|\mathcal{I})}w.
\]
In particular, the one dimensional distribution of the process
$\sqrt{E_\nu(h^2|\mathcal{I})}w$  has a density 
equal to
\[
\frac{1}{2}\frac{1}{\sqrt{2\pi t}}\exp \left(-\frac{x^2}{2t}\right)+
\frac{1}{2}\frac{1}{\sqrt{8\pi t}}\exp
\left(-\frac{x^2}{8t}\right),\;\;x\in\realnos.
\]
\end{example}

In general, for a given $h$ the equation $\mathcal{P}_{T}h=0$ may not be
satisfied. Then the idea is to write $h$ as a sum of two functions, one of which
satisfies the assumptions of Theorem \ref{t:CLT1} while the other is irrelevant
for the convergence to hold. At least a part of the conclusions of
Theorem~\ref{t:CLT2} is given in the following
\begin{theorem}[{\citet[Theorem 3]{tyran}}] \label{t:CLT3} Let $T$ be a non-invertible
measure-preserving transformation
on the probability space $(Y,\B,\nu)$. 
Suppose  $h\in L^2(Y,\B,\nu)$ with $\int h(y)\nu(dy)=0$ is such that
\eqref{concltpo} holds. Then there exists $\tilde{h}\in
L^2(Y,\B,\nu)$ such that $\mathcal{P}_{T}\tilde{h}=0$ and
$\frac{1}{\sqrt{n}}\sum_{j=0}^{n-1}(h-\tilde{h})\circ T^j$ converges
to zero in $L^2(Y,\B,\nu)$ as $n\to\infty$.
\end{theorem}

We will use the following two results for subadditive sequences.

\begin{lemma}[{\citet[Lemma 2.8]{peligradutev}}]\label{l:pu2}
Let $V_n$ be a subadditive sequence of nonnegative numbers. Suppose
that $\sum_{n=1}^{\infty}n^{-3/2}V_n< \infty$. Then
\[
\lim_{m\to\infty}\frac{1}{\sqrt{m}}\sum_{j=0}^{\infty}\frac{V_{m2^j}}{2^{j/2}}=
0.
\]
\end{lemma}

\begin{lemma}\label{l:pu3}
Let $V_n$ be a subadditive sequence of nonnegative numbers. Then for
every integer $r\ge 2$  there exist two positive constants $C_1,C_2$
(depending on $r$) such that
\[
C_1 \sum_{j=0}^{\infty}\frac{V_{r^j}}{r^{j/2}}\le
\sum_{n=1}^{\infty}\frac{V_n}{n^{3/2}}\le C_2
\sum_{j=0}^{\infty}\frac{V_{r^j}}{r^{j/2}}.
\]
\end{lemma}
\begin{proof}
  When $r=2$, the lemma follows from Lemma 2.7 in
  \citet{peligradutev}, the proof of which can be easily extended to the
  case of arbitrary $r>2$.
\end{proof}


\section{Maximal inequality and the proof of
Theorem~\ref{t:CLT2}}\label{s:CLT2}

We start by first stating our key maximal inequality which is analogous to
Proposition~2.3 in \citet{peligradutev}.

\begin{proposition}\label{p:max}
Let $n,q$ be integers such that $2^{q-1}\le n<2^q$. If $T$ is a non-invertible
measure-preserving transformation on the probability space $(Y,\B,\nu)$
and $f\in L^2(Y,\B,\nu)$, then 
\begin{equation}\label{e:max}
    \biggl\lVert \max_{1\le k\le n}\biggl\lvert\sum_{j=0}^{k-1}f\circ
    T^{j}\biggr\rvert\,\biggr\rVert_2\le
    \sqrt{n}\biggl(3\lVert f-U_T\mathcal{P}_{T}f\rVert_2
    +4\sqrt{2}\Delta_q(f)\biggr),
\end{equation}
where
\begin{equation}\label{e:delatar}
\Delta_q(f)=\sum_{j=0}^{q-1}2^{-j/2}\biggl\lVert\sum_{k=1}^{2^j}\mathcal{P}_{T}^k
f\biggr\rVert_2.
\end{equation}
\end{proposition}

In what follows we assume that $T$ is a non-invertible
measure-preserving transformation on the probability space
$(Y,\B,\nu)$.

\begin{proposition}\label{p:thesame}
Let $h\in L^2(Y,\B,\nu)$. Define
\begin{equation}\label{e:samea}
h_m=\frac{1}{\sqrt{m}}\sum_{j=0}^{m-1}h\circ
T^j\qquad\mbox{and}\qquad
w_{k,m}(t)=\frac{1}{\sqrt{k}}\sum_{j=0}^{[kt]-1}h_m\circ T^{mj}
\end{equation}
for $m,k\in\mathbb{N}$ and $t\in[0,1]$. If $m$ is such that the
sequence $ \lVert\max_{1\le l\le k}|w_{k,m}(l/k)|\rVert_2$ is
bounded then
\begin{equation*}
   \lim_{n\to\infty} \bigl\lVert \sup_{0\le t\le 1}|w_{n,1}(t)-w_{[n/m],m}(t)|
    \bigr\rVert_2=0.
\end{equation*}
\end{proposition}
\begin{proof}
Let $k_n=[n/m]$. We have
\[
|w_{n,1}(t)-w_{k_n,m}(t)|\le
\frac{1}{\sqrt{n}}\biggl\lvert\sum_{j=m[k_nt]}^{[nt]-1}h\circ
T^{j}\biggr\rvert +
\biggl(\frac{1}{\sqrt{k_n}}-\frac{\sqrt{m}}{\sqrt{n
}}\biggr)\biggl\lvert\sum_{j=0}^{[k_nt]-1}h_m\circ
T^{mj}\biggr\rvert,
\]
which leads to the estimate
\begin{multline}\label{e:maxw0}
\bigl\lVert \sup_{0\le t\le 1}|w_{n,1}(t)-w_{k_n,m}(t)|
\bigr\rVert_2
       \le \frac{3m}{\sqrt{n}}
    \bigl\lVert \max_{1\le l\le
n}\lvert h\circ T^{l}\rvert
    \bigr\rVert_2 \\
    \quad +\biggl(1-\sqrt{\frac{k_nm}{n}}\biggr)
    \bigl\lVert\max_{1\le l\le k_n}|w_{k_n,m}(l/k_n)|\bigr\rVert_2.
\end{multline}
Since $h\in L^2(Y,\B,\nu)$ we have
\[
\lim_{n\to\infty}\frac{1}{\sqrt{n}}\bigl\lVert \max_{1\le l\le
n}\lvert h\circ T^{l}\rvert
    \bigr\rVert_2=0.
\]
Furthermore, since the sequence $ \lVert\max_{1\le l\le
k}|w_{k,m}(l/k)|\rVert_2$ is bounded by assumption, and
$\lim_{n\to\infty}\bigl(1-\sqrt{ {k_nm}/{n}}\bigr)=0$,
 the second term in the right-hand side of
\eqref{e:maxw0} also tends to $0$.  
\end{proof}

\begin{proof}[Proof of Theorem~\ref{t:CLT2}] From
Theorem~\ref{t:CLT3} it follows that there exists $\tilde{h}\in
L^2(Y,\B,\nu)$ such that $\mathcal{P}_{T}\tilde{h}=0$ and
\begin{equation}\label{e:hht}
    \lim_{n\to\infty}\biggl\lVert\frac{1}{\sqrt{n}}
    \sum_{j=0}^{n-1}(h-\tilde{h})\circ T^j \biggr\rVert_2=0.
\end{equation}
  For each $m\in\mathbb{N}$, define
\[
\tilde{h}_m=\frac{1}{\sqrt{m}}\sum_{j=1}^{m-1}\tilde{h}\circ
T^j\qquad\mbox{and}\qquad
\tilde{w}_{k,m}(t)=\frac{1}{\sqrt{k}}\sum_{j=0}^{[kt]-1}\tilde{h}_m\circ
T^{mj}
\]
for $k\in\mathbb{N}$ and $t\in[0,1]$.  We have
$\mathcal{P}_{T^m}\tilde{h}_m=0$ for all $m$. Thus
Theorem~\ref{t:CLT1} implies
\begin{equation}\label{e:mclt}
\tilde{w}_{k,m}\to^d \sqrt{E_\nu(\tilde{h}_m^2|\mathcal{I}_m)}w
\end{equation}
as $k\to\infty$, where $\mathcal{I}_m$ is the $\sigma$-algebra of
$T^m$-invariant sets. Proposition~\ref{p:max}, applied to $T^m$ and
$\tilde{h}_m$, gives
\[
\bigl \lVert \max_{1\le l\le k}|\tilde{w}_{k,m}(l/k)| \bigr \rVert_2
\le 3 \lVert \tilde{h}_m\rVert_2.
\]
Therefore, by Proposition~\ref{p:thesame}, we obtain
\[
\lim_{n\to\infty}\bigl\lVert \sup_{0\le t\le
1}|\tilde{w}_{n,1}(t)-\tilde{w}_{[n/m],m}(t)|\bigr\rVert_2=0
\]
for all $m\in\mathbb{N}$, which implies, by Theorem~4.1
in~\citet{billingsley68}, that the limit in \eqref{e:mclt} does not
depend on $m$ and is thus equal to
$\sqrt{E_{\nu}(\tilde{h}^2|\mathcal{I})}w$.

To  prove~\eqref{e:clt2}, using Theorem~4.2 in \citet{billingsley68}
 we have to show that
\begin{equation}\label{e:wpsi}
    \lim_{m\to\infty}\limsup_{n\to\infty}
    \bigl\lVert \sup_{0\le t\le 1}|w_n(t)-\tilde{w}_{[n/m],m}(t)|
    \bigr\rVert_2=0.
\end{equation}
Let $h_m$ and  $w_{k,m}$ be defined as in~\eqref{e:samea}. We have
\begin{multline}\label{e:maxw}
\bigl\lVert \sup_{0\le t\le 1}|w_n(t)-\tilde{w}_{[n/m],m}(t)|
\bigr\rVert_2
        \le  \bigl\lVert \sup_{0\le t\le 1}|w_n(t)-w_{[n/m],m}(t)|
\bigr\rVert_2\\
 \quad + \bigl\lVert \sup_{0\le t\le
1}|w_{[n/m],m}(t)-\tilde{w}_{[n/m],m}(t)| \bigr\rVert_2.
\end{multline}
Making use of Proposition~\ref{p:max} with $T^m$ and $h_m$ we obtain
\[
\bigl\lVert \max_{1\le l\le k}\lvert w_{k,m}(l/k)\rvert \bigr
\rVert_2\le 3 \bigl\lVert h_m-U_{T^m}\mathcal{P}_{T^m}
h_m\bigr\rVert_2+4\sqrt{2}\sum_{j=0}^{\infty}2^{-j/2}\biggl\lVert
\sum_{i=1}^{2^j}\mathcal{P}_{T^m}^i h_m\biggr\rVert_2.
\]
However
\[
\mathcal{P}_{T^m}
h_m=\frac{1}{\sqrt{m}}\sum_{j=0}^{m-1}\mathcal{P}_{T^m}U_{T^j}h
=\frac{1}{\sqrt{m}}\sum_{j=1}^{m}\mathcal{P}_{T}^j h,
\]
by~\eqref{fpcond}, and thus
\begin{equation}\label{e:ser}
\sum_{j=0}^{\infty}2^{-j/2}\biggl \lVert
\sum_{i=1}^{2^j}\mathcal{P}_{T^m}^i h_m \biggr \rVert_2=
\frac{1}{\sqrt{m}}\sum_{j=0}^{\infty}2^{-j/2}\biggl \lVert
\sum_{i=1}^{m2^j}\mathcal{P}_{T}^i h \biggr \rVert_2,
\end{equation}
and the series is convergent by Lemma~\ref{l:pu2}, which implies
that the sequence $\bigl\lVert \max_{1\le l\le k}\lvert
w_{k,m}(l/k)\rvert \rVert_2$ is bounded for all $m$. From
Proposition~\ref{p:thesame} it follows that
\[
\lim_{n\to\infty}\bigl\lVert \sup_{0\le t\le
1}|w_n(t)-w_{[n/m],m}(t)| \bigr\rVert_2=0.
\]

We next turn to estimating  the second term in~\eqref{e:maxw}. We
have
\[
\begin{split}
\biggl\lVert \sup_{0\le t\le 1}|w_{k,m}(t)-\tilde{w}_{k,m}(t)|
\biggr\rVert_2 &\le \frac{1}{\sqrt{k}}\biggl\lVert \max_{1\le l\le
k}|\sum_{j=0}^{l-1}(h_m-\tilde{h}_m)\circ T^{mj}| \biggr\rVert_2\\
& \le 3 \bigl\lVert h_m -\tilde{h}_m-U_{T^m}\mathcal{P}_{T^m} (h_m
-\tilde{h}_m)\bigr\rVert_2\\
&\quad +4\sqrt{2}\sum_{j=0}^{\infty}2^{-j/2}\biggl\lVert
\sum_{i=1}^{2^j}\mathcal{P}_{T^m}^i (h_m -\tilde{h}_m)\biggr\rVert_2
\end{split}
\]
by Proposition~\ref{p:max}. Combining this with~\eqref{e:ser} and
the fact that $\mathcal{P}_{T^m}\tilde{h}_m=0$ leads to the estimate
\[
\begin{split}
\biggl\lVert \sup_{0\le t\le 1}|w_{k,m}(t)-\tilde{w}_{k,m}(t)|
\biggr\rVert_2 & \le 3
\frac{1}{\sqrt{m}}\biggl\lVert\sum_{j=0}^{m-1}(h-\tilde{h})\circ T^j
\biggr\rVert_2 +
\frac{1}{\sqrt{m}}\biggl\lVert\sum_{j=1}^m\mathcal{P}_{T^j} h
\biggr\rVert_2\\
&\quad
+\frac{4\sqrt{2}}{\sqrt{m}}\sum_{j=0}^{\infty}2^{-j/2}\biggl\lVert
\sum_{i=1}^{m2^j}\mathcal{P}_{T}^i h\biggr\rVert_2,
\end{split}
\]
which completes the proof of~\eqref{e:wpsi}, because all terms on
the right-hand side tend to $0$ as $m\to\infty$, by~\eqref{e:hht}
and Lemma~\ref{l:pu2}.

\end{proof}

\section{Asymptotically periodic transformations}\label{s:aptrans}

The dynamical properties of what are now known as asymptotically
periodic transformations seem to have first been studied by
\citet{ionescu50}.  These transformations form a perfect example of
the central limit theorem results we have discussed in earlier
sections, and here we consider them in detail.

Let $(X,\A,\mu)$ be a $\sigma$-finite measure space. Let us write
$L^1(\mu)=L^1(X,\A,\mu)$.  The elements of the set
    \[
    D(\mu)=\{f\in L^1(\mu): f\ge 0 \mbox{ and } \int
    f(x)\mu(dx)=1\}
    \]
are called densities. Let $T:X\to X$ be a non-singular
transformation and  $P:L^1(\mu)\to L^1(\mu)$ be the corresponding
Perron-Frobenius operator. 
Then (\citet{almcmbk94}) $(T,\mu)$ is called \emph{asymptotically
periodic} 
if there exists a sequence of densities $g_1,\ldots, g_r$ and a
sequence of bounded linear functionals $\lambda_1,\ldots,\lambda_r$
such that
\begin{equation}\label{d:asp}
\lim_{n\to\infty}\lVert
P^n(f-\sum_{j=1}^r\lambda_j(f)g_j)\rVert_{L^1(\mu)}=0
\end{equation}
for all $f\in D(\mu)$.  The 
densities $g_j$ have disjoint
supports ($g_ig_j=0$ for $i\neq j$) and $Pg_j=g_{\alpha(j)}$, where
$\alpha$ is a permutation of $\{1,\ldots,r\}$.

If $(T,\mu)$ is asymptotically periodic and $r=1$ in \eqref{d:asp} 
then $(T,\mu)$ is called \emph{asymptotically stable} or
\emph{exact} by \citet{almcmbk94}.

Observe that if $(T,\mu)$ is asymptotically periodic then
\[
g_*=\dfrac{1}{r}\sum_{j=1}^rg_j
\]
is an invariant density for $P$, i.e.~$Pg_*=g_*$. The ergodic structure of
asymptotically periodic transformations was studied in~\citet{ishitani91}.


\begin{remark}
Let $\mu(X)<\infty$. Recall that $P$ is a constrictive
Perron-Frobenius  operator if there exists $\delta>0$ and $\kappa<1$
such that for every density $f$
 \[
 \limsup_{n\to\infty}\int_{ A}P^nf(x)\mu(dx)<\kappa
 \]
 for all $A\in\A$ with $\mu(A)\le \delta$.

It is known that if $P$ is a constrictive operator then $(T,\mu)$ is
asymptotically periodic (\citet[Theorem 5.3.1]{almcmbk94}, see also
\citet{komornik87}), and $(T,\mu)$ is ergodic if and only if the
permutation $\{\alpha(1),\ldots,\alpha(r)\}$ of the sequence
$\{1,\ldots,r\}$ is cyclical (\citet[Theorem 5.5.1]{almcmbk94}). In
this case we call $r$ the period of $T$.
\end{remark}


Let  $(T,\mu)$ be asymptotically periodic and let $g_*$ be an
invariant density for $P$. Let $Y=\supp(g_*)=\{x\in X: g_*(x)>0\}$,
$\B=\{A\cap Y:A\in \A\}$, and
\[
\nu(A)=\int_A g_*(x)\mu(dx),\quad A\in\A.
\]
The measure $\nu$ is a probability measure invariant under $T$. In
what follows we write $L^p(\nu)=L^p(Y,\B,\nu)$ for $p=1,2$. The
transfer operator $\mathcal{P}_{T}:L^1(\nu)\to L^1(\nu)$ is given by
\begin{equation}\label{e:topf}
g_*\mathcal{P}_{T}(f)=P(fg_*)\quad \mbox{for}\quad f\in L^1(\nu).
\end{equation}
We now turn to the study of weak convergence of the sequence of
processes
\[
w_{n}(t)=\frac{1}{\sqrt{n}}\sum_{j=0}^{[nt]-1}h\circ T^j,
\]
where $h\in L^2(\nu)$ with $\int h(y)\nu(dy)=0$, by considering
first the ergodic case and then the non-ergodic case.

\subsection{$(T,\mu)$  ergodic and asymptotically
periodic} Let the transformation $(T,\mu)$ be ergodic and
asymptotically periodic with period~$r$. The unique invariant
density of $P$ is given by
$$
g_*=\dfrac{1}{r}\sum_{j=1}^rg_j
$$
and $(T^r,g_j)$ is exact for every $j=1,\ldots,r$. Let
$Y_j=\supp(g_j)$ for $j=1,\ldots, r$. Note that the set
$B_j=\bigcup_{n=0}^{\infty}T^{-nr}(Y_j)$ is (almost) $T^r-$invariant
and $\nu(B_j\setminus Y_j)=0$ for $j=1,\ldots, r$.
Since the $Y_j$ are pairwise disjoint, we have
\[
E_{\nu}(f|\mathcal{I}_r)=\sum_{k=1}^r\dfrac{1}{\nu(Y_k)}\int_{Y_k}f(y)\nu(dy)
1_{Y_k}\quad \mbox{for}\quad f\in L^1(\nu),
\]
where $\mathcal{I}_r$ is the $\sigma$-algebra of $T^r$-invariant
sets. However $\nu(Y_k)=1/r$, and thus
\begin{equation}\label{e:cei}
E_{\nu}(f|\mathcal{I}_r)=r\sum_{k=1}^r\int_{Y_k}f(y)\nu(dy)
1_{Y_k}=\sum_{k=1}^r\int_{Y_k}f(y)g_k(y)\mu(dy)1_{Y_k}.
\end{equation}

\begin{theorem}\label{t:cltap}
Suppose that   $h\in L^2(\nu)$ with $\int h(y)\nu(dy)=0$ is such
that
\begin{equation} \label{concltpor} \sum_{n=1}^\infty
{n^{-\frac{3}{2}}}\biggl\lVert\sum_{k=0}^{n-1}{\mathcal P}_{T}^{rk}
h_r\biggr\rVert_2<\infty,\quad \text{where}\quad
h_r=\dfrac{1}{\sqrt{r}}\sum_{k=0}^{r-1}h\circ T^k.
\end{equation}
Then
\[
w_n\to^d \sigma w,
\]
where $w$ is a standard Brownian motion and $\sigma\ge 0$ is a
constant.

Moreover, if $\sum_{j=1}^\infty\int \lvert
h_r(y)h_r(T^{rj}(y))\rvert\nu(dy)<\infty$ then $\sigma$ is  given by
\begin{equation}
  \sigma^2=r\Bigl(\int_{Y_1} h_r^2(y)\nu(dy)+2\sum_{j=1}^\infty
\int_{Y_1}h_r(y)h_r(T^{rj}(y))\nu(dy)\Bigr).
\end{equation}
\end{theorem}
\begin{proof}
We have $h_r\in L^2(\nu)$ and $\int_{Y}h_r(y)\nu(dy)=0$. Let
\[
w_{k,r}(t)=\frac{1}{\sqrt{k}}\sum_{j=0}^{[kt]-1}h_r\circ T^{rj}\quad
\mbox{for}\quad k\in\mathbb{N},\; t\in[0,1].
\]
We can apply Theorem~\ref{t:CLT2} to deduce that
\[
w_{k,r}\to^d \sqrt{E_{\nu}(\tilde{h}_r^2|\mathcal{I}_r)}w\quad
\mbox{as}\quad k\to\infty,
\]
where $\mathcal{I}_r$ is the $\sigma$-algebra of all $T^r$ invariant
sets and
\begin{equation}\label{e:var}
E_{\nu}(\tilde{h}_r^2|\mathcal{I}_r)=
\lim_{n\to\infty}\frac{1}{n}E_{\nu}\biggl(\bigl(\sum_{j=0}^{n-1}h_r\circ
T^{rj}\bigr)^2|\mathcal{I}_r\biggr).
\end{equation} 
On the other hand, we also have
\[
\sum_{j=0}^{\infty}r^{-j/2}\Bigl\lVert
\sum_{k=1}^{r^j}\mathcal{P}^{rk}h_r\Bigr\rVert_2=
\sum_{j=0}^{\infty}r^{-j/2}\frac{1}{\sqrt{r}}\Bigl\lVert
\sum_{k=1}^{r^{j+1}}\mathcal{P}^{k}h\Bigr\rVert_2=\sum_{j=1}^{\infty}r^{-j/2}\Bigl\lVert
\sum_{k=1}^{r^j}\mathcal{P}^{k}h\Bigr\rVert_2.
\]
Thus the series
\[
\sum_{n=1}^{\infty}n^{-3/2}\Bigl\lVert
\sum_{k=0}^{n-1}\mathcal{P}^{k}h\Bigr\rVert_2
\]
is convergent by Lemma~\ref{l:pu3}. 
From Theorem~\ref{t:CLT2} we conclude that there exists
$\tilde{h}\in L^2(\nu)$ such that
\[
w_n\to^d \lVert \tilde{h}\rVert_2 w
\]
since $T$ is ergodic.  However
\[
\lVert
\tilde{h}\rVert_2=\sqrt{E_{\nu}(\tilde{h}_r^2|\mathcal{I}_r)},
\]
by Proposition~\ref{p:thesame}. Hence
$E_{\nu}(\tilde{h}_r^2|\mathcal{I}_r)$ is a constant and
from~\eqref{e:cei} it follows that for each $k=1,\ldots,r$ the
integral $ \int_{Y_k}\tilde{h}_r^2(y)\nu(dy) $ does not depend on
 $k$. Thus
\[
\sigma^2=\lVert\tilde{h}\rVert_2^2=r\int_{Y_1}\tilde{h}_r^2(y)\nu(dy).
\]
Since $\nu$ is $T^r$-invariant, we have
\[
\begin{split}
\frac{1}{n}\int_{Y_k}
\Bigl(\sum_{j=0}^{n-1}h_r(T^{rj}(y))\Bigr)^2\nu(dy)&=
\int_{Y_k}h_r^2(y)\nu(dy)\\
&\quad +2\frac{1}{n}\sum_{l=1}^{n-1}\sum_{j=1}^{l} \int_{Y_k} h_r(y)
h_r(T^{rj}(y))\nu(dy).
\end{split}
\]
By assumption the sequence $(\sum_{j=1}^n\int_{Y_k}
h_r(y)h_r(T^{rj}(y))\nu(dy))_{n\ge 1}$ is convergent to
$\sum_{j=1}^\infty\int_{Y_k} h_r(y)h_r(T^{rj}(y))\nu(dy)$ which
completes the proof when combined with~\eqref{e:var} and
\eqref{e:cei}.
\end{proof}


\subsection{${{(T,\mu)}}$ asymptotically periodic but not
necessarily ergodic} Now let us consider $(T,\mu)$ asymptotically
periodic but not ergodic, so that the permutation $\alpha$ is not
cyclical and we can represent it as a product of permutation cycles.
Thus we can rephrase the definition of asymptotic periodicity as
follows.

Let there exist a sequence of densities
\begin{equation}
g_{1,1},\ldots, g_{1,r_1}, \ldots, g_{l,1},\ldots,g_{l,r_l}
\end{equation} and a sequence
of bounded linear functionals $\lambda_{1,1},\ldots,
\lambda_{1,r_1}, \ldots, \lambda_{l,1},\ldots,\lambda_{l,r_l}$ such
that
\begin{equation}\label{e:asp1}
\lim_{n\to\infty}\lVert
P^n(f-\sum_{i=1}^l\sum_{j=1}^{r_i}\lambda_{i,j}(f)g_{i,j})\rVert_{L^1(\mu)}=0
\quad\text{for all}\quad f\in L^1(\mu),
\end{equation}
where the densities $g_{i,j}$ have mutually disjoint supports and
for each $i$, $Pg_{i,j}=g_{i,j+1}$ for $1\le j\le r_i-1$,
$Pg_{i,r_i}=g_{i,1}$. 
Then
\[
g_i^*=\frac{1}{r_i}\sum_{j=1}^{r_i}g_{i,j}
\]
is an invariant density for $P$  and $(T,g_i^*)$ is ergodic for
every $i=1,\ldots, l$. Let $g_*$ be a convex combination of $g_i^*$,
i.e.
\[
g_*=\sum_{i=1}^l \alpha_i g_i^*
\]
where $\alpha_i\ge 0$ and $\sum_{i=1}^{l}\alpha_i=1$. 
For simplicity, assume that $\alpha_i>0$.

Let 
$Y_i=\supp(g_i^*)$ and $Y_{i,j}=\supp(g_{i,j})$, $j=1,\ldots,r_i$, $i=1,\ldots,l$.
If $\mathcal{I}$ is the $\sigma$-algebra of all $T$-invariant sets, then
\[
E_{\nu}(f|\mathcal{I})=\sum_{i=1}^l\frac{1}{\nu(Y_i)}\int_{Y_i}
f(y)\nu(dy)1_{Y_i}=\sum_{i=1}^l\int_{Y_i}
f(y)g_i^*(y)\mu(dy)1_{Y_i}.
\]
Now, if $\mathcal{I}_r$ is the $\sigma$-algebra of all
$T^r$-invariant sets with $r=\prod_{i=1}^lr_i$, then
\[
E_{\nu}(f|\mathcal{I}_r)=\sum_{i=1}^l\frac{r_i}{\nu(Y_i)}\sum_{j=1}^{r_i}\int_{Y_{i,j}}
f(y)\nu(dy)1_{Y_{i,j}}
\]
for $f\in L^1(\nu)$, which leads to
\[
E_{\nu}(f|\mathcal{I}_r)=\sum_{i=1}^l\sum_{j=1}^{r_i}\int_{Y_{i,j}}
f(y)g_{i,j}(y)\mu(dy)1_{Y_{i,j}}.
\]
Using similar arguments as in the proof of  Theorem \ref{t:cltap} we
obtain
\begin{theorem}\label{t:cltap2}
Suppose that   $h\in L^2(\nu)$ with $\int h(y)\nu(dy)=0$ is such
that condition~\eqref{concltpor} holds.
Then
\[
w_n\to^d \eta w,
\]
where $w$ is a standard Brownian motion and $\eta\ge 0$ is a random
variable independent of $w$.

Moreover, if $\sum_{j=1}^\infty\int \lvert
h_r(y)h_r(T^{rj}(y))\rvert\nu(dy)<\infty$ then $\eta$ is given by
\begin{equation*}
  \eta=\sum_{i=1}^l\frac{r_i}{\nu(Y_i)}\Bigl(\int_{Y_{i,1}} h_r^2(y)\nu(dy)+2\sum_{j=1}^\infty
\int_{Y_{i,1}}h_r(y)h_r(T^{rj}(y))\nu(dy)\Bigr)1_{Y_i}.
\end{equation*}
\end{theorem}

%
%

%

\begin{remark}\label{r:cond}
Observe that condition~\eqref{concltpor} holds if
\begin{equation*}
\sum_{n=1}^\infty \frac{\lVert{\mathcal P}_{T}^{rn}
h_r\rVert_2}{\sqrt{n}}<\infty.
\end{equation*}
The operator ${\mathcal P}_{T}$ is a contraction on $L^\infty(\nu)$.
Therefore \begin{equation*} \lVert\mathcal{P}_{T}^n f\rVert_2\le
\lVert
f\lVert_\infty^{1/2}\lVert\mathcal{P}_{T}^nf\rVert_1^{1/2}
\quad\text{for}\quad f\in L^\infty(\nu), \;n\ge 1,
\end{equation*}
which allows us to easily check condition~\eqref{concltpor} for
specific examples of transformations $T$.

It also should be noted that, by \eqref{e:topf}, we have
\[
\lVert\mathcal{P}_{T}^nf\rVert_1=\lVert{
P}^n(fg_*)\rVert_{L^1(\mu)}\quad \text{for}\quad f\in
L^1(\nu),\;n\ge 1.
\]
\end{remark}

\subsection{Piecewise monotonic transformations}

Let $X$ be a totally ordered, order complete set (usually $X$ is a compact interval in $\mathbb{R}$). 
Let $\B$ be the $\sigma$-algebra of Borel subsets of $X$ and let
$\mu$ be a  probability measure on $X$. Recall that a function
$f:X\to\mathbb{R}$ is said to be of {\it bounded variation} if
\[
\mathrm{var}(f)=\sup\sum_{i=1}^{n}|f(x_{i-1})-f(x_i)|<\infty,
\]
where the supremum is taken over all finite ordered sequences,
$(x_j)$ with $x_j\in X$.
%
The bounded variation norm is given by
\[
||f||_{BV}=\lVert f\rVert_{L^1(\mu)}+\mathrm{var}(f)
\]
and it makes $BV = \{f : X\to\mathbb{R} : \mathrm{var}(f) <\infty\}$
into a Banach space.

Let $T : V\to X$ be a continuous map, $V\subset X$ be open and dense
with $\mu(V) = 1$. We call $(T,\mu)$ a \emph{piecewise uniformly
expanding map} if:
\begin{enumerate}
\item There exists a countable family $\Z$ of closed intervals
with disjoint interiors such that $V\subset\bigcup_{Z\in\Z}Z$ and for any $Z\in\Z$
the set $Z\cap (X\setminus V)$ consists exactly of the endpoints of $Z$.
\item For any $Z\in \Z$, $T_{|Z\cap V}$
admits an extension to a homeomorphism from $Z$ to some interval.
\item There exists a function $g : X\to [0,\infty)$, with bounded variation,
$g_{|X\setminus V }= 0$ such that the Perron-Frobenius operator $P :
L^1(\mu) \to L^1(\mu)$ is of the form
\[Pf (x) = \sum_{z\in T^{-1}(x)} g(z)f (z).\]
\item $T$ is expanding: $\sup_{x\in V} g(x)<1$.
\end{enumerate}

The following result is due to \citet{rychlik}
\begin{theorem}
If $(T,\mu)$ is a piecewise uniformly expanding map then 
it satisfies \eqref{e:asp1} with $g_{i,j}\in BV$. Moreover, there
exist constants $C>0$ and $\theta\in(0,1)$ such that for every
function $f$ of bounded variation and all $n\ge 1$
\[
\lVert {P}^{rn}f-Q(f)\rVert_{L^1(\mu)}\le C \theta^n||f ||_{BV},
\]
where $r=\prod_{i=1}^l r_i$  and
\[
Q(f)=\sum_{i=1}^l\sum_{j=1}^{r_i} \int_{Y_{i,j}}f(x)\mu(dx) g_{i,j}.
\]

%
\end{theorem}

This result and Remark~\ref{r:cond} imply
\begin{corollary}
Let $(T,\mu)$ be a piecewise uniformly expanding map and $\nu$ an
invariant measure which is absolutely continuous with respect to
measure $\mu$. If $h$ is a function of bounded
variation with $E_{\nu}(h|\mathcal{I})=0$ then condition~\eqref{concltpor} holds.  
\end{corollary}

\begin{remark}
AFU-maps (Uniformly expanding maps satisfying Adler's condition with
a Finite image condition, which are interval maps with a finite
number of indifferent fixed points) studied in \citet{zweimuller98},
are asymptotically periodic when they have an absolutely continuous
invariant probability measure. However, the decay of the $L^1$ norm
may not be exponential. For H\"older continuous functions $h$ one
might use the results of \citet{young99} to obtain bounds on this
norm and then apply our results.
\end{remark}

\subsection{Calculation of variance for the family
of tent maps using Theorem~\ref{t:cltap}}\label{s:tent}

Let  $T$ be  the generalized tent map on $[-1,1]$ defined by
\begin{equation}
T_a(x) = a - 1 -a |x| \qquad \mbox{for} \quad x \in \left [ -1, 1
\right ], \label{eqn:tent2}
 \end{equation} where $a\in (1,2]$.
The Perron-Frobenius operator $P:L^1(\mu)\to L^1(\mu)$ is given by
\begin{equation}\label{e:FPtent}
  Pf(x)=\dfrac{1}{a}\left(f\left(\psi_a^{-}(x)\right)
  +f\left(\psi_a^{+}(x)\right)\right)1_{[-1,a-1]}(x),
\end{equation}
where $\psi_a^{-}$ and $\psi_a^{+}$ are the inverse branches of
$T_a$
\begin{equation}
\psi_a^{-}(x)=\dfrac{x+1-a}{a},\qquad
\psi_a^{+}(x)=-\dfrac{x+1-a}{a}
\end{equation}
and $\mu$ is the normalized Lebesgue measure on $[-1,1]$.

 \citet{ito79a} have  shown that the tent  map
Equation \ref{eqn:tent2}  is ergodic, thus possessing a unique
invariant density $g_a$. 
\citet{provatas91a} have proved the asymptotic periodicity of
(\ref{eqn:tent2}) with period $r = 2^m$ for
\begin{equation*}
2^{1/2^{ {m+1}} } < a \le 2^{1/2^{ {m}} }\qquad \mbox{for} \qquad
m=0,1,2,\cdots. 
\end{equation*} Thus, for example, $(T,\mu)$ has period $1$ for $2^{1/2} < a
\leq 2$, period $2$ for $2^{1/4} < a \leq 2^{1/2}$, period $4$ for
$2^{1/8} < a \leq 2^{1/4}$, etc.

Let $Y=\supp(g_a)$ and $\nu_a(dy)=g_a(y)\mu(dy)$. For all $1<a\le 2$
we have $T_a(A)=A$ with $A=[T_a^2(0),T_a(0)]$ and $g_a(x)=0$ for
$x\in[-1,1]\setminus A$. If $\sqrt{2}<a\le 2$ then $g_a$ is strictly
positive in $A$, thus $Y=A$ in this case. For $a\le \sqrt{2}$ we
have $Y\subset A$. The transfer operator $\mathcal{P}_a\colon
L^1(\nu_a)\to L^1(\nu_a)$ is given by
\begin{equation*}\label{e:transoper}
\mathcal{P}_af=\dfrac{P(fg_a)}{g_a}\quad\mbox{for}\quad f\in
L^1(\nu_a),
\end{equation*}
where $P$ is the Perron-Frobenius operator \eqref{e:FPtent}.

If $h$ is a function of bounded variation on $[-1,1]$ with
$\int_{-1}^{1}h(y)\nu_a(dy)=0$ and
\[
w_{n}(t)=\frac{1}{\sqrt{n}}\sum_{j=0}^{[nt]-1}h\circ T_a^j,
\]
then there exists a constant $\sigma(h)\ge 0$  such that
\begin{equation*}
 w_n\to^d \sigma(h) w,
\end{equation*}
where $w$ is a standard Brownian motion. In particular, we are going
to study $\sigma(h)$ for the specific example of $h=h_a$ for $a\in
(1,2]$ , where
\begin{equation*}
h_{a}(y)=y-\mathfrak{m}_a, \;
y\in[-1,1], 
 \quad \mbox{and}\quad \mathfrak{m}_a=\int_{[-1,1]}yg_a(y)\,dy.
\end{equation*}

\begin{proposition}\label{p:tent}
Let $m\ge 1$ and $2^{1/2^{m+1}} < a \le 2^{1/2^m}$. Then
\begin{equation}\label{e:sigmaai}
\sigma(h_a)=\frac{\sigma(h_{a^{2^m}})a(a-1)}{\sqrt{2^m}a^{2^m}(a^{2^m}-1)}\prod_{k=0}^{m-1}(a^{2^k}-1)^2,
\end{equation}
where
\begin{equation}\label{e:sigmaha2m}
\begin{split}
 \sigma(h_{a^{2^m}})^2&=
2 \int h_{a^{2^m}}(y) f_{a^{2^m}}(y) \nu_{a^{2^m}}(dy)-\int h_{a^{2^m}}^2(y)\nu_{a^{2^m}}(dy)\\
\quad&\mbox{and}\quad
f_{a^{2^m}}=\sum_{n=0}^{\infty}\mathcal{P}_{a^{2^m}}^nh_{a^{2^m}}.
\end{split}
\end{equation}
\end{proposition}

In general, an explicit representation  for~\eqref{e:sigmaha2m} is
not known.  Hence, before turning to a proof of Proposition
\ref{p:tent}, we first give the simplest example in which
$\sigma(h_{a^{2^m}})^2$ can be calculated exactly.

\begin{example}
For $a=2$ the invariant density for the transformation $T_a$ is
$g_2=\frac{1}{2}1_{[-1,1]}$ and the transfer operator
$\mathcal{P}_2\colon L^1(\nu_2)\to L^1(\nu_2)$ has the same form as
$P$ in \eqref{e:FPtent}  \[ \mathcal{P}_2f=\frac{1}{2}\bigl(f\circ
\psi_2^{-}+f\circ \psi_2^+\bigr).
\]
Since $\int_{-1}^{1} y dy=0$, we have $h_2(y)=y$. We also have
$\mathcal{P}_2h_2=0$. Thus \[\sigma(h_2)^2=\frac{1}{2}\int_{-1}^{1}
y^2 dy=1/3\] and Proposition~\ref{p:tent} gives $\sigma(h_a)$ for
$a=2^{1/2^m}$, $m\ge 1$.
%
\end{example}

We now summarize  some properties of the tent map \cite{yoshida83},
which will
allow us to prove Proposition~\ref{p:tent}. 
Let $I_0=[x^*(a),x^*(a)(1+\frac{2}{a})]$ and $I_1=[-x^*(a),x^*(a)]$,
where  $x^*(a)$ is the fixed point of $T_a$ other than $-1$, i.e.
$$x^*(a)=\dfrac{a-1}{a+1}.$$ Define transformations
$\phi_{ia}:I_i\to[-1,1]$ by
$$\phi_{1a}(x)=-\dfrac{1}{x^*(a)}x\quad \mbox{and}\quad
\phi_{0a}(x)=\dfrac{a}{x^*(a)}x-a-1.$$ We have
\begin{equation}\label{e:phiia}
\phi_{1a}^{-1}(x)=-x^*(a)x\quad \mbox{and}\quad
\phi_{0a}^{-1}(x)=\dfrac{x^*(a)}{a}\bigl(x+a+1\bigr).
\end{equation}
Then for $1<a\le \sqrt{2}$ the map $T_a^2:I_i\to I_i$ is conjugate
to $T_{a^2}:[-1,1]\to[-1,1]$
\begin{equation}\label{e:cona2}
T_{a^2}=\phi_{ia}\circ T_a^2\circ \phi_{ia}^{-1}
\end{equation} and the invariant density of $T_a$ is given by
\begin{equation}\label{e:inga}
  g_a(y)=\dfrac{1}{2x^*(a)}\left(a
  g_{a^2}(\phi_{0a}(y))1_{I_0}(y)+g_{a^2}(\phi_{1a}(y))1_{I_1}(y))\right).
\end{equation}

\begin{lemma}\label{l:ch2}
If $a\in (1,\sqrt{2}]$ then
\begin{equation}\label{e:pmean}
\mathfrak{m}_a=\dfrac{a-1}{2a}-\dfrac{(a-1)x^*(a)}{2a}\mathfrak{m}_{a^2}
\end{equation}
and
\begin{equation}\label{e:haha2}
(h_{a}+h_{a}\circ T_a)\circ
\phi_{0a}^{-1}=\dfrac{(1-a)x^*(a)}{a}h_{a^2}
\end{equation}
\end{lemma}
\begin{proof}
Equation \eqref{e:pmean} follows from \eqref{e:inga} and
\eqref{e:phiia}, while \eqref{e:haha2} is a direct consequence of
the definition of $\phi_{0a}^{-1}$, the fact that $I_0\subset
[0,1]$, and  \eqref{e:pmean}.
\end{proof}


 Let $m\ge 1$. For
$2^{1/2^{m+1}}<a\le 2^{1/2^{m}}$ there exist $2^m$ disjoint
intervals in which $g_a$ is strictly positive and they are defined
by
\begin{equation*}
  Y_j^m=\Phi_{jm}^{-1}([T^2_{a^{2^m}}(0),T_{a^{2^m}}(0)]),
\end{equation*}
where 
\[
\Phi_{jm}=\phi_{i_m a^{2^{m-1}}}\circ \phi_{i_{m-1}
a^{2^{m-2}}}\circ\ldots \phi_{i_2 a^2}\circ \phi_{i_1 a}
\]
 and $j=1+i_1+2i_2+\ldots+2^{m-1}i_m$, $i_k=0,1$,
$k=1,\ldots,m$. We have $T_a(Y_j^m)=Y_{j+1}^m$ for $1\le j\le 2^m-1$
and $T_a(Y_{2^m}^m)=Y_1^m$. In particular,  we have
\begin{equation}\label{e:Y1m}
Y_1^{m+1}=\phi_{0a}^{-1}(Y_1^{m})\quad \text{for}\quad m\ge 0,
\end{equation} 
where  $Y_1^0=[T^2_{a^{2}}(0),T_{a^{2}}(0)]$. 


\begin{lemma}\label{l:indcor}
Define
\begin{equation}\label{e:hra}
h_{r,a}=\dfrac{1}{\sqrt{r}}\sum_{k=0}^{r-1}h_a\circ T_a^k
\quad\text{for}\quad r\ge 1, \;a\in(1,2].
\end{equation}
Let $m\ge 0$ and $r=2^{m}$. If $2^{1/4r}<a\le 2^{1/2r}$ then
\begin{equation}\label{e:chind}
\begin{split}
 \int_{Y_1^{m+1}}h_{2r,a}(y)h_{2r,a}&(T_a^{2rn}(y))\nu_a(dy)=\\
&\frac{(1-a)^2x^*(a)^2}{2^2a^2}\int_{Y_1^{m}}h_{r,a^2}(y)h_{r,a^2}(T_{a^2}^{rn}(y))\nu_{a^2}(dy)
\end{split}
\end{equation}
for all $n\ge 0$.
\end{lemma}
\begin{proof}
First observe that
\begin{equation}\label{e:ch2}
h_{2r,a}=\dfrac{1}{\sqrt{r}}\sum_{k=0}^{r-1}h_{2,a}\circ T_a^{2k}.
\end{equation}
Let $n\ge 0$.  Since $\phi_{0a}^{-1}(\phi_{0a}(y))=y$ for
$y\in[-1,1]$, a change of variables using \eqref{e:Y1m} and
\eqref{e:inga} gives
\begin{multline}\label{e:chva}
\int_{Y_1^{m+1}}h_{2r,a}(y)h_{2r,a}(T_a^{2rn}(y))\nu_a(dy)=
\\
 \frac{1}{2}\int_{Y_1^m}h_{2r,a}(\phi_{0a}^{-1}(y))h_{2r,a}(T_a^{2rn}(\phi_{0a}^{-1}(y)))\nu_{a^2}
 (dy).
\end{multline}
We have $T_a^{2k}\circ\phi_{0a}^{-1}=\phi_{0a}^{-1}\circ
T_{a^2}^{k}$ for all $k\ge 0$ by \eqref{e:cona2}. Thus
$T_a^{2rn}\circ\phi_{0a}^{-1}=\phi_{0a}^{-1}\circ T_{a^2}^{rn}$ and
from \eqref{e:ch2}  it follows that
\[
h_{2r,a}\circ\phi_{0a}^{-1}
=\dfrac{1}{\sqrt{r}}\sum_{k=0}^{r-1}h_{2,a}\circ\phi_{0a}^{-1}\circ
T_{a^2}^{k}.
\]
By Lemma~\ref{l:ch2} we obtain
\[
h_{2,a}\circ\phi_{0a}^{-1}=\frac{(1-a)x^*(a)}{\sqrt{2}a}h_{a^2}.
\]
Hence
\begin{equation*}
h_{2r,a}\circ\phi_{0a}^{-1}=\frac{(1-a)x^*(a)}{\sqrt{2}a}h_{r,a^2},
\end{equation*}
which, when substituted into equation \eqref{e:chva}, completes the
proof.
\end{proof}

\begin{proof}[Proof of Proposition \ref{p:tent}]
First, we show that if  $m\ge 1$ and $2^{1/2^{m+1}}<a\le
2^{1/2^{m}}$ then
\begin{equation}\label{e:sigmaaii}
\sigma(h_a)=\frac{\sigma(h_{a^{2^m}})}{\sqrt{2^m}a^{2^m-1}}
\prod_{k=0}^{m-1}x^*(a^{2^k})(a^{2^k}-1).
\end{equation}
Let $m\ge 1$ and $2^{1/2^{m+1}}<a\le 2^{1/2^{m}}$. Since the
transformation $T_a$ is asymptotically periodic with period $2^{m}$,
Theorem~\ref{t:cltap} gives
\[
\sigma(h_a)^2=2^{m}\Bigl(\int_{Y_1^m}
h_{2^{m},a}^2(y)\nu_a(dy)+2\sum_{j=1}^\infty
\int_{Y_1^{m}}h_{2^{m},a}(y)h_{2^{m},a}(T_a^{2^mj}(y))\nu_a(dy)\Bigr).
\]
We have $a^2\in (2^{1/2^m},2^{1/2^{m-1}}]$ and the transformation
$T_{a^2}$ is asymptotically periodic with period $r=2^{m-1}$. From
\eqref{e:chind} with $r=2^{m-1}$ and Theorem~\ref{t:cltap} it
follows that
\begin{equation*}
\sigma(h_a)^2=\frac{(a-1)^2x^*(a)^2}{2a^2}\sigma(h_{a^2})^2.
\end{equation*}
Thus equation \eqref{e:sigmaaii} follows immediately by an induction
argument on $m$. Finally, we have for each $k=0,\ldots,m-1$
\[
x^*(a^{2^k})(a^{2^k}-1)=\frac{a^{2^k}-1}{a^{2^k}+1}(a^{2^k}-1)
=\frac{(a^{2^k}-1)^3}{a^{2^{k+1}}-1}
\]
and equation \eqref{e:sigmaai} holds. Since $a^{2^m}>\sqrt{2}$ the
function $f_{a^{2^m}}$ is well defined and
\[
\int h_{a^{2^m}}(y) f_{a^{2^m}}(y)
\nu_{a^{2^m}}(dy)=\sum_{n=0}^\infty \int h_{a^{2^m}}(y)
h_{a^{2^m}}(T_{a^{2^m}}^n(y))\nu_{a^{2^m}}(dy),
\]
which completes the proof.
\end{proof}

\section*{Acknowledgments} This work was supported by the Natural
Sciences and  Engineering Research Council (NSERC grant OGP-0036920,
Canada) and the Mathematics of Information Technology and Complex
Systems (MITACS Canada). This research was carried out while MCM was
visiting University of Silesia, and MT-K was visiting McGill
University.

\appendix

\section{Proof of the maximal inequality}\label{a:app}

\begin{proof}[Proof of Proposition \ref{p:max}]
We will prove \eqref{e:max} inductively. If $n=1$ and $q=1$ we have
\[
\lVert f\rVert_2\le \lVert f-U_T\mathcal{P}_{T}f\rVert_2 + \lVert
 U_T\mathcal{P}_{T}f\rVert_2=\lVert f-U_T\mathcal{P}_{T}f\rVert_2+\Delta_{1}(f)
\]
by the invariance of $\nu$ under $T$.  Now assume that \eqref{e:max} holds for all
$n< 2^{q-1}$. Fix $n$, $2^{q-1}\le n<2^q$. By the triangle inequality
\begin{equation}\label{e:trian}
\max_{1\le k\le n}\biggl\lvert\sum_{j=0}^{k-1}f\circ
    T^{j}\biggr\rvert\le
    \max_{1\le k\le n}\biggl\lvert\sum_{j=0}^{k-1}(f-U_T\mathcal{P}_{T}f)\circ
    T^{j}\biggr\rvert
    +\max_{1\le k\le n}\biggl\lvert\sum_{j=0}^{k-1}U_T\mathcal{P}_{T}f\circ
    T^{j}\biggr\rvert.
    \end{equation}
We first show that
 \begin{equation}\label{e:msec}
 \biggl\lVert\max_{1\le k\le n}\biggl\lvert
\sum_{j=0}^{k-1}(f-U_T\mathcal{P}_{T}f)\circ
    T^{j}\biggr\rvert \biggr\rVert_2 \le 3\sqrt{n}\lVert
    f-U_T\mathcal{P}_{T}f\rVert_2.
    \end{equation}
Observe  that
\[
\begin{split}
 \max_{1\le k\le n}\biggl\lvert\sum_{j=0}^{k-1}(f-U_T\mathcal{P}_{T}f)\circ
    T^{j}\biggr\rvert&\le \biggl\lvert\sum_{j=0}^{n-1}(f-U_T\mathcal{P}_{T}f)\circ
    T^{j}\biggr\rvert\\
    &\quad + \max_{1\le k\le n}\biggl\lvert\sum_{j=1}^{k}(f-U_T\mathcal{P}_{T}f)\circ
    T^{n-j}\biggr\rvert.
    \end{split}
\]
Since $\mathcal{P}_{T}(f-U_T\mathcal{P}_{T}f)=0$, we see that
\[
\biggl\lVert\sum_{j=0}^{n-1}(f-U_T\mathcal{P}_{T}f)\circ
    T^{j}\biggr\rVert_2=\sqrt{n}\bigl\lVert f-U_T\mathcal{P}_{T}f
    \bigr\rVert_2.
\]
 For every $n$ the family $\{\sum_{j=1}^k(f-U_T\mathcal{P}_{T}f)\circ
    T^{n-j}: 1\le k\le n\}$ is
a martingale with respect to $\{T^{-n+k}(\B):1\le k\le n\}$.  Thus by the Doob
maximal inequality
     \begin{equation*}
    \begin{split}
\biggl\lVert\max_{1\le k\le n}\biggl\lvert
\sum_{j=1}^{k}(f-U_T\mathcal{P}_{T}f)\circ
    T^{n-j}\biggr\rvert \biggr\rVert_2 &\le 2
    \biggl\lVert \sum_{j=1}^{n}(f-U_T\mathcal{P}_{T}f)\circ
    T^{n-j}\biggr\rVert_2\\
    & =2\sqrt{n}\lVert f-U_T\mathcal{P}_{T}f\rVert_2,
    \end{split}
    \end{equation*}
which completes the proof of  \eqref{e:msec}.

Now consider the second term on the right hand side of~\eqref{e:trian}. Writing
$n=2m$ or $n=2m+1$ yields
\begin{equation}\label{e:msec1}
\max_{1\le k\le n}\biggl\lvert\sum_{j=0}^{k-1}U_T\mathcal{P}_{T}f\circ
    T^{j}\biggr\rvert\le \max_{1\le l\le m}\biggl\lvert{\sum_{j=0}^{l-1}f_1\circ
    T^{2j}}\biggr\rvert +\max_{0\le l\le m}\biggl\lvert U_T\mathcal{P}_{T}f\circ
    T^{2l}\biggr\rvert,
\end{equation}
where $f_1=U_{T^2}\mathcal{P}_{T}f+U_T\mathcal{P}_{T}f$. To estimate the norm of
the second term in the right hand side of \eqref{e:msec1}, observe that
\[
\max_{0\le l\le m}\lvert U_T\mathcal{P}_{T}f\circ
    T^{2l}\rvert^2\le \sum_{l=0}^{m}\lvert U_T\mathcal{P}_{T}f\circ
    T^{2l}\rvert^2,
\]
which leads to
\begin{equation}\label{e:msec2}
   \biggl\lVert \max_{0\le l\le m}\lvert U_T\mathcal{P}_{T}f\circ
    T^{2l}\rvert \biggr\rVert_2\le\sqrt{m+1}\lVert \mathcal{P}_{T}
    f\rVert_2,
\end{equation}
since $\nu$ is invariant under $T$. Further, since $m<2^{q-1}$, the measure $\nu$
is invariant under $T^2$, and $f_1\in L^2(Y,\B,\nu)$, we can use the induction
hypothesis. We thus obtain
\[
\biggl\lVert\max_{1\le l\le m}\biggl\lvert\sum_{j=0}^{l-1}f_1\circ
    T^{2j}\biggr\rvert\biggr\rVert_2\le
    \sqrt{m}\biggl(3\lVert f_1-U_{T^2}\mathcal{P}_{T^2}f_1\rVert_2
    +4\sqrt{2}\Delta_{q-1}(f_1)\biggr).
\]
We have $f_1-U_{T^2}\mathcal{P}_{T^2}f_1
=U_T\mathcal{P}_{T}f-U_{T^2}\mathcal{P}_{T^2}f$, by \eqref{fpcond}, which implies
\[
\lVert f_1-U_{T^2}\mathcal{P}_{T^2}f_1\rVert_2\le
\lVert\mathcal{P}_{T}f\rVert_2+\lVert\mathcal{P}_{T^2}f\rVert_2\le
2\lVert\mathcal{P}_{T}f\rVert_2, 
\]
since $\mathcal{P}_{T}$ is a contraction. We also have
\[
\begin{split}
\Delta_{q-1}(f_1)&=\sum_{j=0}^{q-2}2^{-j/2}\biggl\lVert\sum_{k=1}^{2^j}
\mathcal{P}_{T^2}^k f_1\biggr\rVert_2
=\sum_{j=0}^{q-2}2^{-j/2}\biggl\lVert\sum_{k=1}^{2^j}
\mathcal{P}_{T}^{2k}  f_1\biggr\rVert_2\\
&=\sum_{j=0}^{q-2}2^{-j/2}\biggl\lVert\sum_{k=1}^{2^j} \mathcal{P}_{T}^{2k}
(U_{T^2}\mathcal{P}_{T}f+U_T\mathcal{P}_{T}f)\biggr\rVert_2\\
&=\sum_{j=0}^{q-2}2^{-j/2}\biggl\lVert\sum_{k=1}^{2^{j+1}} \mathcal{P}_{T}^k
f\biggr\rVert_2=\sqrt{2}\biggl(\Delta_{q}(f)-\lVert
\mathcal{P}_{T}f\rVert_2\biggr).
\end{split}
\]
Therefore
\[
\biggl\lVert\max_{1\le l\le m} \biggl \lvert {\sum_{j=0}^{l-1}f_1\circ
    T^{2j}} \biggr \rvert  \biggr\rVert_2\le \sqrt{m}\bigl(8\Delta_{q}(f)-2\lVert
\mathcal{P}_{T}f\rVert_2\bigr),
\]
which combined with \eqref{e:trian}
through \eqref{e:msec2} and the fact that $\sqrt{m+1}\le \sqrt{2m}\le\sqrt{n}$
leads to
\[
\begin{split}
\biggl \lVert \max_{1\le k\le n} \biggl \lvert \sum_{j=1}^{k}f\circ
    T^{n-j} \biggr \rvert \biggr \rVert_2 & \le 3\sqrt{n}\lVert f-U_T\mathcal{P}_{T}f\rVert_2
    + \sqrt{m+1} \lVert \mathcal{P}_{T}f\rVert_2\\
    &\quad +\sqrt{2m}\bigl(4\sqrt{2}\Delta_{q}(f)-\sqrt{2}\lVert
\mathcal{P}_{T}f\rVert_2\bigr)\\
&\le \sqrt{n} \bigl(3\lVert
f-U_T\mathcal{P}_{T}f\rVert_2+4\sqrt{2}\Delta_{q}(f)\bigr).\qedhere
    \end{split}
\]
\end{proof}

\section{The limiting random variable $\eta$}\label{b:app}

Finally, we give a series expansion of $E_\nu(\tilde{h}^2|\mathcal{I})$ in
Theorem~\ref{t:CLT2} in terms of $h$ and iterates of $T$.

\begin{proposition}
Suppose  $h\in L^2(Y,\B,\nu)$ with $\int h(y)\nu(dy)=0$ is such that
\begin{equation}\label{e:sums}
  \sum_{j=0}^{\infty}2^{-j/2}\lVert\sum_{k=1}^{2^j} \mathcal{P}_{T}^k
  h\rVert_2<\infty.
\end{equation}
Then the following limit exists in $L^1$
\begin{equation}\label{e:lim}
 \lim_{n\to\infty}\frac{E_{\nu}(S_n^2|\mathcal{I})}{n}=
  E_{\nu}(h^2|\mathcal{I})+\sum_{j=0}^{\infty}
  \frac{E_{\nu}(S_{2^j}S_{2^j}\circ T^{2^j}|\mathcal{I})}{2^j},
\end{equation}
where $\mathcal{I}$ is the $\sigma$-algebra of all $T$-invariant sets and
$S_{n}=\sum_{j=0}^{n-1}h\circ T^j$, $n\in\mathbb{N}$.

Moreover, if $\tilde{h}\in L^2(Y,\B,\nu)$ is such that
$\mathcal{P}_{T}\tilde{h}=0$ and
$\bigl\lVert\frac{1}{\sqrt{n}}\sum_{j=0}^{n-1}(h-\tilde{h})\circ
T^j\bigr\rVert_2\to 0$ as $n\to\infty$ then
\begin{equation}\label{e:limt}
E_{\nu}(\tilde{h}^2|\mathcal{I})=\lim_{n\to\infty}\frac{E_{\nu}(S_n^2|\mathcal{I})}{n}.
\end{equation}
\end{proposition}

\begin{proof}
We first prove that the series in the right-hand side of~\eqref{e:lim} is
convergent in $L^1(Y,\B,\nu)$. Since $\mathcal{I}\subset T^{-2^j}(\B)$ for all
$j$, we see that
\[
E_{\nu}(S_{2^j}S_{2^j}\circ
T^{2^j}|\mathcal{I})=E_{\nu}(E_{\nu}(S_{2^j}S_{2^j}\circ
T^{2^j}|T^{-2^j}(\B))|\mathcal{I}).
\]
As $S_{2^j}\circ T^{2^j}$ is $T^{-2^j}(\B)$-measurable and integrable we have
\[E_{\nu}(S_{2^j}S_{2^j}\circ T^{2^j}|T^{-2^j}(\B))=S_{2^j}\circ
T^{2^j}E_{\nu}(S_{2^j}|T^{-2^j}(\B)).\] 
However,
$E_{\nu}(S_{2^j}|T^{-2^j}(\B))=U_{T}^{2^{j}}\mathcal{P}_{T}^{2^{j}}S_{2^j}$ from
~\eqref{fpcond}. Consequently,
\begin{equation}\label{e:s2j}
E_{\nu}(S_{2^j}S_{2^j}\circ
T^{2^j}|\mathcal{I})=E_\nu(S_{2^j}\sum_{k=1}^{2^j}\mathcal{P}_{T}^k
h|\mathcal{I}).
\end{equation}
Since the conditional expectation operator is a contraction in $L^1$, we have
\[
\lVert E_{\nu}(S_{2^j}S_{2^j}\circ T^{2^j}|\mathcal{I})\rVert_1 \le
\lVert S_{2^j}\sum_{k=1}^{2^j}\mathcal{P}_{T}^k h\rVert_1,
\]
which, by the Cauchy-Schwartz inequality, leads to 
\[
\lVert E_{\nu}(S_{2^j}S_{2^j}\circ T^{2^j}|\mathcal{I})\rVert_1\le\lVert
S_{2^j}\rVert_2
 \lVert\sum_{k=1}^{2^j}\mathcal{P}_{T}^k h\rVert_2.
\]
Since $\lVert S_{2^j}\rVert_2\le \lVert \max_{1\le l\le 2^j}\lvert S_{l}\rvert
\rVert_2$, the sequence $\lVert S_{2^j}\rVert_2/2^{j/2}$ is bounded
by~\eqref{e:sums}, Lemma~\ref{l:pu3}, and Proposition~\ref{p:max}. Hence
\[
\sum_{j=0}^{\infty}\frac{\lVert S_{2^j}\rVert_2
 \lVert\sum_{k=1}^{2^j}\mathcal{P}_{T}^k h\rVert_2}{2^j}\le C\sum_{j=0}^{\infty}\frac{
 \lVert\sum_{k=1}^{2^j}\mathcal{P}_{T}^k h\rVert_2}{2^{j/2}}<\infty,
\]
which  proves the convergence in $L^1$ of the series in~\eqref{e:lim}.

We now prove the equality in \eqref{e:lim}. Since
\[
S_{2^m}^2=\bigl(S_{2^{m-1}}+S_{2^{m-1}}\circ
T^{2^{m-1}}\bigr)^2=S_{2^{m-1}}^2+S_{2^{m-1}}^2\circ
T^{2^{m-1}}+2S_{2^{m-1}}S_{2^{m-1}}\circ T^{2^{m-1}},
\]
we obtain
\[
E_{\nu}(S_{2^m}^2|\mathcal{I})=2E_{\nu}(S_{2^{m-1}}^2|\mathcal{I})
+2E_{\nu}(S_{2^{m-1}}S_{2^{m-1}}\circ T^{2^{m-1}}|\mathcal{I}),
\]
which leads to
\[
\frac{E_{\nu}(S_{2^m}^2|\mathcal{I})}{2^m
}=E_{\nu}(h^2|\mathcal{I})+\sum_{j=0}^{m-1}\frac{E_{\nu}(S_{2^j}S_{2^j}\circ
T^{2^j}|\mathcal{I})}{2^j}.
\]
Thus the limit in the left-hand side of~\eqref{e:lim} exists for the subsequence
$n=2^m$ and the equality holds. An analysis similar to that in the proof of
Proposition~2.1 in \citet{peligradutev} shows that the whole sequence is
convergent, which completes the proof of
\eqref{e:lim}. 

We now turn to the proof of~\eqref{e:limt}.  Let $\tilde{h}$ be such that
$\mathcal{P}_{T}\tilde{h}=0$. Define $\tilde{S}_n=\sum_{j=0}^{n-1}\tilde{h}\circ
 T^j$. Substituting $\tilde{h}$
into~\eqref{e:sums} and~\eqref{e:s2j} gives
\[
 E_{\nu}(\tilde{h}^2|\mathcal{I})=\lim_{n\to\infty}\frac{E_{\nu}(\tilde{S}_n^2|\mathcal{I})}{n}.
\]
We have
\[
\biggl\lVert\frac{E_{\nu}(\tilde{S}_n^2|\mathcal{I})}{n}
-\frac{E_{\nu}(S_n^2|\mathcal{I})}{n}\biggr\rVert_1
\le\biggl\lVert\frac{\tilde{S}_n^2}{n}-\frac{S_n^2}{n}\biggr\rVert_1\le
\biggl\lVert\frac{\tilde{S}_n}{\sqrt{n}}-\frac{S_n}{\sqrt{n}}\biggr\rVert_2
\biggl\lVert\frac{\tilde{S}_n}{\sqrt{n}}+\frac{S_n}{\sqrt{n}}\biggr\rVert_2
\]
by the H\"older inequality, which implies \eqref{e:limt} when combined with 
the equality $\bigl\lVert\sum_{j=0}^{n-1}\tilde{h}\circ
T^j\bigr\rVert_2=\sqrt{n}\lVert \tilde{h}\rVert_2$, and the assumption
$\bigl\lVert\frac{1}{\sqrt{n}}\sum_{j=0}^{n-1}(h-\tilde{h})\circ
T^j\bigr\rVert_2\to 0$ as $n\to\infty$.
\end{proof}

\bibliographystyle{apalike}
\bibliography{c:/my-stuff/PAPERS/BIB-FILES/zpf}

\begin{thebibliography}{}

\bibitem[Billingsley, 1968]{billingsley68}
Billingsley, P. (1968).
\newblock {\em Convergence of probability measures}.
\newblock John Wiley \& Sons Inc., New York.

\bibitem[Billingsley, 1995]{billingsley95}
Billingsley, P. (1995).
\newblock {\em Probability and measure}.
\newblock Wiley Series in Probability and Mathematical Statistics. John Wiley
  \& Sons Inc., New York.

\bibitem[Boyarsky and Scarowsky, 1979]{boyarsky79}
Boyarsky, A. and Scarowsky, M. (1979).
\newblock On a class of transformations which have unique absolutely continuous
  invariant measures.
\newblock {\em Trans. Amer. Math. Soc.}, 255:243--262.

\bibitem[Conze and Le~Borgne, 2001]{conze01}
Conze, J.-P. and Le~Borgne, S. (2001).
\newblock M\'ethode de martingales et flot g\'eod\'esique sur une surface de
  courbure constante n\'egative.
\newblock {\em Ergodic Theory Dynam. Systems}, 21:421--441.

\bibitem[Denker, 1989]{denker89}
Denker, M. (1989).
\newblock The central limit theorem for dynamical systems.
\newblock In {\em Dynamical systems and ergodic theory (Warsaw, 1986)},
  volume~23 of {\em Banach Center Publ.}, pages 33--62. PWN, Warsaw.

\bibitem[Eagleson, 1975]{eagleson75}
Eagleson, G.~K. (1975).
\newblock Martingale convergence to mixtures of infinitely divisible laws.
\newblock {\em Ann. Probab.}, 3(3):557--562.

\bibitem[Gordin, 1969]{gordin69}
Gordin, M.~I. (1969).
\newblock The central limit theorem for stationary processes.
\newblock {\em Dokl. Akad. Nauk SSSR}, 188:739--741.

\bibitem[Gou{\"e}zel, 2004]{gouezel04b}
Gou{\"e}zel, S. (2004).
\newblock Central limit theorem and stable laws for intermittent maps.
\newblock {\em Probab. Theory Related Fields}, 128(1):82--122.

\bibitem[Inoue and Ishitani, 1991]{ishitani91}
Inoue, T. and Ishitani, H. (1991).
\newblock Asymptotic periodicity of densities and ergodic properties for
  nonsingular systems.
\newblock {\em Hiroshima Math. J.}, 21(3):597--620.

\bibitem[Ionescu~Tulcea and Marinescu, 1950]{ionescu50}
Ionescu~Tulcea, C.~T. and Marinescu, G. (1950).
\newblock Th\'eorie ergodique pour des classes d'op\'erations non
  compl\`etement continues.
\newblock {\em Ann. of Math. (2)}, 52:140--147.

\bibitem[Ito et~al., 1979]{ito79a}
Ito, S., Tanaka, S., and Nakada, H. (1979).
\newblock On unimodal linear transformations and chaos. {I}.
\newblock {\em Tokyo J. Math.}, 2(2):221--239.

\bibitem[Jab{\l}o{\'n}ski and Malczak, 1983]{jab83}
Jab{\l}o{\'n}ski, M. and Malczak, J. (1983).
\newblock A central limit theorem for piecewise convex mappings of the unit
  interval.
\newblock {\em T\^ohoku Math. J. (2)}, 35(2):173--180.

\bibitem[Keller, 1980]{keller80}
Keller, G. (1980).
\newblock Un th\'eor\`eme de la limite centrale pour une classe de
  transformations monotones par morceaux.
\newblock {\em C. R. Acad. Sci. Paris S\'er. A-B}, 291(2):A155--A158.

\bibitem[Komorn{\'{\i}}k and Lasota, 1987]{komornik87}
Komorn{\'{\i}}k, J. and Lasota, A. (1987).
\newblock Asymptotic decomposition of {M}arkov operators.
\newblock {\em Bull. Polish Acad. Sci. Math.}, 35(5-6):321--327.

\bibitem[Lasota and Mackey, 1994]{almcmbk94}
Lasota, A. and Mackey, M.~C. (1994).
\newblock {\em Chaos, fractals, and noise}, volume~97 of {\em Applied
  Mathematical Sciences}.
\newblock Springer-Verlag, New York.

\bibitem[Liverani, 1996]{liverani}
Liverani, C. (1996).
\newblock Central limit theorem for deterministic systems.
\newblock In {\em International Conference on Dynamical Systems (Montevideo,
  1995)}, volume 362 of {\em Pitman Res. Notes Math. Ser.}, pages 56--75.
  Longman, Harlow.

\bibitem[Mackey and Tyran-Kami{\'n}ska, 2006]{mackeytyran}
Mackey, M.~C. and Tyran-Kami{\'n}ska, M. (2006).
\newblock Deterministic {B}rownian motion: {T}he effects of perturbing a
  dynamical system by a chaotic semi-dynamical system.
\newblock {\em Phys. Rep.}, 422(5):167--222.

\bibitem[Maxwell and Woodroofe, 2000]{maxwell00}
Maxwell, M. and Woodroofe, M. (2000).
\newblock Central limit theorems for additive functionals of {M}arkov chains.
\newblock {\em Ann. Probab.}, 28(2):713--724.

\bibitem[Melbourne and Nicol, 2004]{melbourne04}
Melbourne, I. and Nicol, M. (2004).
\newblock Statistical properties of endomorphisms and compact group extensions.
\newblock {\em J. London Math. Soc. (2)}, 70(2):427--446.

\bibitem[Melbourne and T{\"o}r{\"o}k, 2002]{melbourne}
Melbourne, I. and T{\"o}r{\"o}k, A. (2002).
\newblock Central limit theorems and invariance principles for time-one maps of
  hyperbolic flows.
\newblock {\em Comm. Math. Phys.}, 229(1):57--71.

\bibitem[Merlev{\`e}de et~al., 2006]{merlevedeetal06}
Merlev{\`e}de, F., Peligrad, M., and Utev, S. (2006).
\newblock Recent advances in invariance principles for stationary sequences.
\newblock {\em Probab. Surv.}, 3:1--36 (electronic).

\bibitem[Peligrad and Utev, 2005]{peligradutev}
Peligrad, M. and Utev, S. (2005).
\newblock A new maximal inequality and invariance principle for stationary
  sequences.
\newblock {\em Ann. Probab.}, 33(2):798--815.

\bibitem[Peligrad et~al., 2006]{peligradutevwu}
Peligrad, M., Utev, S., and Wu, W.~B. (2006).
\newblock A maximal ${L}_p$-inequality for stationary sequences and its
  applications.
\newblock {\em Proc. Amer. Math. Soc.}, in press.

\bibitem[Provatas and Mackey, 1991]{provatas91a}
Provatas, N. and Mackey, M.~C. (1991).
\newblock Asymptotic periodicity and banded chaos.
\newblock {\em Phys. D}, 53(2-4):295--318.

\bibitem[Rousseau-Egele, 1983]{rousseau-egele}
Rousseau-Egele, J. (1983).
\newblock Un th\'eor\`eme de la limite locale pour une classe de
  transformations dilatantes et monotones par morceaux.
\newblock {\em Ann. Probab.}, 11(3):772--788.

\bibitem[Rychlik, 1983]{rychlik}
Rychlik, M. (1983).
\newblock Bounded variation and invariant measures.
\newblock {\em Studia Math.}, 76(1):69--80.

\bibitem[Tyran-Kami{\'n}ska, 2005]{tyran}
Tyran-Kami{\'n}ska, M. (2005).
\newblock An invariance principle for maps with polynomial decay of
  correlations.
\newblock {\em Comm. Math. Phys.}, 260(1):1--15.

\bibitem[Voln{\'y}, 1987a]{volny87b}
Voln{\'y}, D. (1987a).
\newblock A nonergodic version of {G}ordin's {CLT} for integrable stationary
  processes.
\newblock {\em Comment. Math. Univ. Carolin.}, 28(3):413--419.

\bibitem[Voln{\'y}, 1987b]{volny87}
Voln{\'y}, D. (1987b).
\newblock On the invariance principle and functional law of iterated logarithm
  for nonergodic processes.
\newblock {\em Yokohama Math. J.}, 35(1-2):137--141.

\bibitem[Voln{\'y}, 1989]{volny89}
Voln{\'y}, D. (1989).
\newblock On nonergodic versions of limit theorems.
\newblock {\em Apl. Mat.}, 34(5):351--363.

\bibitem[Voln{\'y}, 1993]{volny93}
Voln{\'y}, D. (1993).
\newblock Approximating martingales and the central limit theorem for strictly
  stationary processes.
\newblock {\em Stochastic Process. Appl.}, 44(1):41--74.

\bibitem[Wong, 1979]{wong79}
Wong, S. (1979).
\newblock A central limit theorem for piecewise monotonic mappings of the unit
  interval.
\newblock {\em Ann. Probab.}, 7(3):500--514.

\bibitem[Yoshida et~al., 1983]{yoshida83}
Yoshida, T., Mori, H., and Shigematsu, H. (1983).
\newblock Analytic study of chaos of the tent map: band structures, power
  spectra, and critical behaviors.
\newblock {\em J. Statist. Phys.}, 31(2):279--308.

\bibitem[Young, 1999]{young99}
Young, L.-S. (1999).
\newblock Recurrence times and rates of mixing.
\newblock {\em Israel J. Math.}, 110:153--188.

\bibitem[Zweim{\"u}ller, 1998]{zweimuller98}
Zweim{\"u}ller, R. (1998).
\newblock Ergodic structure and invariant densities of non-{M}arkovian interval
  maps with indifferent fixed points.
\newblock {\em Nonlinearity}, 11(5):1263--1276.

\end{thebibliography}
\end{document}